\documentclass[11pt, a4paper]{article}
\usepackage[cp1251]{inputenc}
\usepackage{amsmath} \usepackage{euscript} \usepackage{marvosym}
\usepackage[dvipsnames]{xcolor}
\usepackage{mymatrx6}
\usepackage{graphicx}
\usepackage{longtable}
\usepackage{mathrsfs}
\usepackage{amssymb, amscd}
\usepackage{comment}
\def\mymatrix{\MyMatrixwithdelims..}
\Linestrue\Autonumfalse
\begin{document}

\newcounter{bnomer} \newcounter{snomer}
\newcounter{bsnomer}
\setcounter{bnomer}{0}
\renewcommand{\thesnomer}{\thebnomer.\arabic{snomer}}
\renewcommand{\thebsnomer}{\thebnomer.\arabic{bsnomer}}
\renewcommand{\refname}{\begin{center}\large{\textbf{References}}\end{center}}

\setcounter{MaxMatrixCols}{14}

\newcommand\restr[2]{{
  \left.\kern-\nulldelimiterspace 
  #1 
  \right|_{#2} 
}}

\newcommand{\sect}[1]{%
\setcounter{snomer}{0}\setcounter{bsnomer}{0}
\refstepcounter{bnomer}
\par\bigskip\begin{center}\large{\textbf{\arabic{bnomer}. {#1}}}\end{center}}
\newcommand{\sst}[1]{%
\refstepcounter{bsnomer}
\par\bigskip\textbf{\arabic{bnomer}.\arabic{bsnomer}. {#1}}\par}
\newcommand{\defi}[1]{%
\refstepcounter{snomer}
\par\medskip\textbf{Definition \arabic{bnomer}.\arabic{snomer}. }{#1}\par\medskip}
\newcommand{\theo}[2]{%
\refstepcounter{snomer}
\par\textbf{Theorem \arabic{bnomer}.\arabic{snomer}. }{#2} {\emph{#1}}\hspace{\fill}$\square$\par}
\newcommand{\mtheop}[2]{%
\refstepcounter{snomer}
\par\textbf{Theorem \arabic{bnomer}.\arabic{snomer}. }{\emph{#1}}
\par\textsc{Proof}. {#2}\hspace{\fill}$\square$\par}
\newcommand{\mcorop}[2]{%
\refstepcounter{snomer}
\par\textbf{Corollary \arabic{bnomer}.\arabic{snomer}. }{\emph{#1}}
\par\textsc{Proof}. {#2}\hspace{\fill}$\square$\par}
\newcommand{\mtheo}[1]{%
\refstepcounter{snomer}
\par\medskip\textbf{Theorem \arabic{bnomer}.\arabic{snomer}. }{\emph{#1}}\par\medskip}
\newcommand{\theobn}[1]{%
\par\medskip\textbf{Theorem. }{\emph{#1}}\par\medskip}
\newcommand{\theoc}[2]{%
\refstepcounter{snomer}
\par\medskip\textbf{Theorem \arabic{bnomer}.\arabic{snomer}. }{#1} {\emph{#2}}\par\medskip}
\newcommand{\mlemm}[1]{%
\refstepcounter{snomer}
\par\medskip\textbf{Lemma \arabic{bnomer}.\arabic{snomer}. }{\emph{#1}}\par\medskip}
\newcommand{\mprop}[1]{%
\refstepcounter{snomer}
\par\medskip\textbf{Proposition \arabic{bnomer}.\arabic{snomer}. }{\emph{#1}}\par\medskip}
\newcommand{\theobp}[2]{%
\refstepcounter{snomer}
\par\textbf{Theorem \arabic{bnomer}.\arabic{snomer}. }{#2} {\emph{#1}}\par}
\newcommand{\theop}[2]{%
\refstepcounter{snomer}
\par\textbf{Theorem \arabic{bnomer}.\arabic{snomer}. }{\emph{#1}}
\par\textsc{Proof}. {#2}\hspace{\fill}$\square$\par}
\newcommand{\theosp}[2]{%
\refstepcounter{snomer}
\par\textbf{Theorem \arabic{bnomer}.\arabic{snomer}. }{\emph{#1}}
\par\textsc{Sketch of the proof}. {#2}\hspace{\fill}$\square$\par}
\newcommand{\exam}[1]{%
\refstepcounter{snomer}
\par\medskip\textbf{Example \arabic{bnomer}.\arabic{snomer}. }{#1}\par\medskip}
\newcommand{\deno}[1]{%
\refstepcounter{snomer}
\par\textbf{Notation \arabic{bnomer}.\arabic{snomer}. }{#1}\par}
\newcommand{\lemm}[1]{%
\refstepcounter{snomer}
\par\textbf{Lemma \arabic{bnomer}.\arabic{snomer}. }{\emph{#1}}\hspace{\fill}$\square$\par}
\newcommand{\lemmp}[2]{%
\refstepcounter{snomer}
\par\medskip\textbf{Lemma \arabic{bnomer}.\arabic{snomer}. }{\emph{#1}}
\par\textsc{Proof}. {#2}\hspace{\fill}$\square$\par\medskip}
\newcommand{\coro}[1]{%
\refstepcounter{snomer}
\par\textbf{Corollary \arabic{bnomer}.\arabic{snomer}. }{#1}\hspace{\fill}$\square$\par}
\newcommand{\mcoro}[1]{%
\refstepcounter{snomer}
\par\textbf{Corollary \arabic{bnomer}.\arabic{snomer}. }{\emph{#1}}\par\medskip}
\newcommand{\corop}[2]{%
\refstepcounter{snomer}
\par\textbf{Corollary \arabic{bnomer}.\arabic{snomer}. }\emph{#1}
\par\textsc{Proof}. {#2}\hspace{\fill}$\square$\par}
\newcommand{\nota}[1]{%
\refstepcounter{snomer}
\par\medskip\textbf{Remark \arabic{bnomer}.\arabic{snomer}. }{#1}\par\medskip}
\newcommand{\propp}[2]{%
\refstepcounter{snomer}
\par\medskip\textbf{Proposition \arabic{bnomer}.\arabic{snomer}. }{\emph{#1}}
\par\textsc{Proof}. {#2}\hspace{\fill}$\square$\par\medskip}
\newcommand{\hypo}[1]{%
\refstepcounter{snomer}
\par\medskip\textbf{Conjecture \arabic{bnomer}.\arabic{snomer}. }{\emph{#1}}\par\medskip}
\newcommand{\prop}[1]{%
\refstepcounter{snomer}
\par\textbf{Proposition \arabic{bnomer}.\arabic{snomer}. }{\emph{#1}}\hspace{\fill}$\square$\par}

\newcommand{\proof}[2]{%
\par\medskip\textsc{Proof{#1}}. \hspace{-0.2cm}{#2}\hspace{\fill}$\square$\par\medskip}

\makeatletter
\def\iddots{\mathinner{\mkern1mu\raise\p@
\vbox{\kern7\p@\hbox{.}}\mkern2mu
\raise4\p@\hbox{.}\mkern2mu\raise7\p@\hbox{.}\mkern1mu}}
\makeatother

\newcommand{\okr}[2]{%
\refstepcounter{snomer}
\par\medskip\textbf{{#1} \arabic{bnomer}.\arabic{snomer}. }{\emph{#2}}\par\medskip}

\newcommand{\Ind}[3]{%
\mathrm{Ind}_{#1}^{#2}{#3}}
\newcommand{\Res}[3]{%
\mathrm{Res}_{#1}^{#2}{#3}}
\newcommand{\epsi}{\varepsilon}
\newcommand{\tri}{\triangleleft}
\newcommand{\Supp}[1]{%
\mathrm{Supp}(#1)}
\newcommand{\SSu}[1]{%
\mathrm{SingSupp}(#1)}

\newcommand{\gee}{\geqslant}
\newcommand{\reg}{\mathrm{reg}}
\newcommand{\Dyn}{\mathrm{Dyn}}
\newcommand{\Ann}{\mathrm{Ann}\,}
\newcommand{\Cent}[1]{\mathbin\mathrm{Cent}({#1})}
\newcommand{\PCent}[1]{\mathbin\mathrm{PCent}({#1})}
\newcommand{\Irr}[1]{\mathbin\mathrm{Irr}({#1})}
\newcommand{\Exp}[1]{\mathbin\mathrm{Exp}({#1})}
\newcommand{\empr}[2]{[-{#1},{#1}]\times[-{#2},{#2}]}
\newcommand{\sreg}{\mathrm{sreg}}
\newcommand{\ilm}{\varinjlim}
\newcommand{\wdth}{\mathrm{wd}}
\newcommand{\plm}{\varprojlim}
\newcommand{\codim}{\mathrm{codim}\,}
\newcommand{\GKdim}{\mathrm{GKdim}\,}
\newcommand{\chara}{\mathrm{char}\,}
\newcommand{\rk}{\mathrm{rk}\,}
\newcommand{\chr}{\mathrm{ch}\,}
\newcommand{\Ker}{\mathrm{Ker}\,}
\newcommand{\id}{\mathrm{id}}
\newcommand{\Ad}{\mathrm{Ad}}
\newcommand{\Gh}{\mathrm{Gh}}
\newcommand{\col}{\mathrm{col}}
\newcommand{\row}{\mathrm{row}}
\newcommand{\high}{\mathrm{high}}
\newcommand{\low}{\mathrm{low}}
\newcommand{\pho}{\hphantom{\quad}\vphantom{\mid}}
\newcommand{\fho}[1]{\vphantom{\mid}\setbox0\hbox{00}\hbox to \wd0{\hss\ensuremath{#1}\hss}}
\newcommand{\wt}{\widetilde}
\newcommand{\wh}{\widehat}
\newcommand{\ad}[1]{\mathrm{ad}_{#1}}
\newcommand{\tr}{\mathrm{tr}\,}
\newcommand{\GL}{\mathrm{GL}}
\newcommand{\SL}{\mathrm{SL}}
\newcommand{\SO}{\mathrm{SO}}
\newcommand{\Or}{\mathrm{O}}
\newcommand{\Sp}{\mathrm{Sp}}
\newcommand{\Sa}{\mathrm{S}}
\newcommand{\Ua}{\mathrm{U}}
\newcommand{\Andre}{\mathrm{Andre}}
\newcommand{\Aord}{\mathrm{Aord}}
\newcommand{\Mat}{\mathrm{Mat}}
\newcommand{\Stab}{\mathrm{Stab}}
\newcommand{\htt}{\mathfrak{h}}
\newcommand{\spt}{\mathfrak{sp}}
\newcommand{\slt}{\mathfrak{sl}}
\newcommand{\sot}{\mathfrak{so}}

\newcommand{\vfi}{\varphi}
\newcommand{\aad}{\mathrm{ad}}
\newcommand{\vpi}{\varpi}
\newcommand{\teta}{\vartheta}
\newcommand{\Bfi}{\Phi}
\newcommand{\Fp}{\mathbb{F}}
\newcommand{\Rp}{\mathbb{R}}
\newcommand{\Zp}{\mathbb{Z}}
\newcommand{\Cp}{\mathbb{C}}
\newcommand{\Ap}{\mathbb{A}}
\newcommand{\Pp}{\mathbb{P}}
\newcommand{\Kp}{\mathbb{K}}
\newcommand{\Np}{\mathbb{N}}
\newcommand{\ut}{\mathfrak{u}}
\newcommand{\at}{\mathfrak{a}}
\newcommand{\glt}{\mathfrak{gl}}
\newcommand{\hei}{\mathfrak{hei}}
\newcommand{\nt}{\mathfrak{n}}
\newcommand{\kt}{\mathfrak{k}}
\newcommand{\mt}{\mathfrak{m}}
\newcommand{\rt}{\mathfrak{r}}
\newcommand{\rad}{\mathfrak{rad}}
\newcommand{\bt}{\mathfrak{b}}
\newcommand{\unt}{\underline{\mathfrak{n}}}
\newcommand{\gt}{\mathfrak{g}}
\newcommand{\vt}{\mathfrak{v}}
\newcommand{\pt}{\mathfrak{p}}
\newcommand{\Xt}{\mathfrak{X}}
\newcommand{\Po}{\mathcal{P}}
\newcommand{\PV}{\mathcal{PV}}
\newcommand{\Uo}{\EuScript{U}}
\newcommand{\Fo}{\EuScript{F}}
\newcommand{\Do}{\EuScript{D}}
\newcommand{\Eo}{\EuScript{E}}
\newcommand{\Jo}{\EuScript{J}}
\newcommand{\Iu}{\mathcal{I}}
\newcommand{\Mo}{\mathcal{M}}
\newcommand{\Nu}{\mathcal{N}}
\newcommand{\Ro}{\mathcal{R}}
\newcommand{\Co}{\mathcal{C}}
\newcommand{\Ko}{\mathcal{K}}
\newcommand{\So}{\mathcal{S}}
\newcommand{\Lo}{\mathcal{L}}
\newcommand{\Ou}{\mathcal{O}}
\newcommand{\Uu}{\mathcal{U}}
\newcommand{\Tu}{\mathcal{T}}
\newcommand{\Au}{\mathcal{A}}
\newcommand{\Vu}{\mathcal{V}}
\newcommand{\Du}{\mathcal{D}}
\newcommand{\Bu}{\mathcal{B}}
\newcommand{\Sy}{\mathcal{Z}}
\newcommand{\Sb}{\mathcal{F}}
\newcommand{\Gr}{\mathcal{G}}
\newcommand{\Xu}{\mathcal{X}}
\newcommand{\Op}{\mathbb{O}}
\newcommand{\chv}{\mathrm{chv}}
\newcommand{\rtc}[1]{C_{#1}^{\mathrm{red}}}

\author{Mikhail Ignatev\and Mikhail Venchakov}
\date{}
\title{Orbits of maximal and submaximal dimension\\ for Sylow $p$-subgroups of finite classical groups}\maketitle
\begin{abstract} Let $U$ be a Sylow $p$-subgroup in a classical group over a finite field with $q$ elements of characteristic $p$ large enough. The coadjoint orbits of the group $U$ play the key role in the description of irreducible complex characters of $U$. In the paper, we provide a classification of such orbits of pre-maximal dimension for symplectic groups and orbits of maximal dimension for orthogonal groups. As a corollary, we compute the number of all orbits mentioned above. It turned out that each of these numbers is a polynomial in $q - 1$ with integer non-negative coefficients, which agrees with the Isaacs' conjecture.


\end{abstract}

\sect{Introduction}

\let\thefootnote\relax\footnote{The research was supported by RSF (project No. 25--21--00219), \texttt{https://rscf.ru/en/project/25-21-00219/}.}

Let $U$ be a unipotent algebraic group over a finite field $\Fp_q$ of sufficiently large characteristic $p$. The main tool in representation theory of $U$ is the orbit method created in 1962 by A.A. Kirillov, see~\cite{Kirillov62}, \cite{Kirillov04}, \cite{Kazhdan77}. The key idea of the orbit method says that the irreducible representations of $U$ are in one-to-one correspondence with the coadjoint orbits of this group. Namely, the group $U$ acts on its Lie algebra $\ut$ via the adjoint action: $$g\in U,~x\in\ut\mapsto gxg^{-1}\in U;$$ the dual action of $U$ on the dual space $\ut^*$ is called coadjoint. It turns out that there is a natural bijection between the set $\Irr{U}$ of all irreducible complex characters of $U$ and the set $\ut^*/U$ of coadjoint orbits.

Let $U$ be a maximal unipotent subgroup (or, equivalently, a Sylow $p$-subgroup) in a simple classical group over $\Fp_q$. A complete description of all coadjoint orbits for the group $U$ is a wild problem, so a natural question is how to describe certain important classes of orbits and the corresponding irreducible characters. For the case $A_{n-1}$, orbits of maximal possible dimension (we call such orbits regular) were classified in the first Kirillov's work on the orbit method \cite{Kirillov62}. All of them are associated with so-called orthogonal rook placements.

For other classical root systems, so-called orbits associated with the Kostant cascades have maximal possible dimension, see \cite{Kostant12} and \cite{Kostant13}, but not all the orbits of maximal dimension have such a form. For $A_{n-1}$, orbits of submaximal dimension (we call such orbits subregular) were classified by A.N. Panov \cite{IgnatevPanov09}. Such orbits also correspond to orthogonal rook placements (modulo adding simple root covectors to the canonical forms on them). A classification of orbits of maximal dimension in type $C_n$ follows from C. Andre and A. Neto's paper \cite{AndreNeto06}, see also \cite{Venchakov}. A description of orbits of maximal possible dimension for the root systems $B_n$ and $D_n$, as well as of orbits of submaximal dimension for all classical types $B_n$, $C_n$ and $D_n$, was still unknown. 

\newpage
The main goal of this paper is to fill the gap particularly. 
The structure of the paper is as follows.
In the next section we introduce a necessary notations and definitions. The main results are formulated in the next sections. More precisely, in the third section we provide an explicit classification of the subregular orbits (i.e., the orbits of submaximal dimension) for the type $C_n$, prove Theorem~\ref{theo_orb_max_dim} and calculate their number as a corollary. The fourth section is dedicated to the description of the maximal orbits for groups of type $B_n$ and $D_n$, see Theorem~\ref{theo_orb_max_dim_BD}. The main technical tool which was used is the method of C-patterns and C-quatterns, invented in \cite{GoodwinMoschRoehrle16} by S.M. Goodwin, P. Mosch and G. R\"ohrle, and then applied to description of orbits of certain special classes in \cite{IgnatevPetukhov25} by the first author and A. Petukhov and \cite{Venchakov} by the second author. Notice also that this method can not be used directly to obtain a classification of subregular orbits of groups of type $B_n$ and $D_n$ (see Remark~\ref{nota1}).

A longstanding G. Higman's conjecture \cite{Higman60} states that the number of conjugacy classes for~$U=U_n(q)$, the group of all upper-triangular matrices over $\Fp_q$ with 1's on the diagonal (which is a Sylow $p$-subgroup in $\SL_n(\Fp_q)$), is a polynomial in~$q$. It is easy to check that the number of conjugacy classes for $U$ coincides with the number of coadjoint $U$-orbits, or, equivalently, with the number of irreducible characters of $U$. Of course, it is interesting to consider an analogue of Higman's conjecture for an arbitrary $U$, not only for $U_n(q)$; we will denote the number of irreducible characters of the group $U$ by $O(q)$. Fourteen years after Higman, G. Lehrer conjectured \cite{Lehrer74} that, for $U=U_n(q)$, even the number $O_e(q)$ of characters of degree $q^e$ is polynomial in $q$. In 2007, M. Isaacs \cite{Isaacs07} made a stronger conjecture that, for $U=U_n(q)$, $O_e(q)$ is in fact polynomial in $q-1$ with nonnegative integer coefficients. In past twenty five years, a significant progress in studying these conjectures has been made for $U=U_n(q)$, see, e.g., \cite{VeraLopezArregi03}, \cite{Marjoram97'}, \cite{Marjoram97}, \cite{IgnatevPanov09}, \cite{Marjoram99}, \cite{Loukaki11}, \cite{Le10}, \cite{Marberg11}.

For Sylow subgroups $U$ in other Chevalley groups, the situation is as follows. In 1999, Marjoram computed $O_{\mu(\Phi)}(q)$ for orthogonal group (i.e., for $\Phi=B_n$ or $D_n$), where $q^{\mu(\Phi)}$ is the maximal possible degree of an irreducible character of $U$, see \cite{Marjoram99}. For the symplectic case (i.e., for $\Phi=C_n$), the formula for $O_{\mu(\Phi)}(q)$ follows from the results of C.A.M. Andr\`e and A.-M. Neto on so-called supercharacters published in 2006 \cite{AndreNeto06}. In 2016, S.M. Goodwin, P. Mosch and G. R\"ohrle calculated $O_e(q)$ and proved Isaacs's conjecture for all possible $e$ and all finite Chevalley groups $G(q)$ of rank $\leq 8$, except $E_8$. Our classification of regular orbits for $B_n$ and $D_n$ gives another proof for the result of Marjoram (see Corollary \ref{coro:Isaacs_BD}), while the classification of subregular orbits for $C_n$ gives a new evidence for the Isaacs' conjecture in this case (see Corollary \ref{coro:Isaacs_C}).

We thank Alexey Petukhov for very useful discussions.


\sect{Main definitions} \label{main_definitions}
To begin with we present some basic facts about root systems of simple algebras. We will denote root systems of type $B_n$, $D_n$ or $C_n$ by $\Phi$. As usual, we will identify them with the following subsets of $\Rp^n$:
$$
\begin{array}{ll}
&B_n = \{\pm\epsi_i\pm\epsi_j, 1\leq i < j\leq
n\}\cup\{\pm\epsi_i, 1\leq i\leq n\},\\
&D_n = \{\pm\epsi_i\pm\epsi_j, 1\leq i < j\leq n\},\\
&C_n = \{\pm\epsi_i\pm\epsi_j, 1\leq i < j\leq
n\}\cup\{\pm2\epsi_i, 1\leq i\leq n\},
\end{array}
$$
where $\{\epsi_i\}_{i=1}^n$ is the standard basis in $\Rp^n$. Pick the set of simple roots $\Delta=\Delta(\Phi)$ as in
\cite{Bourbaki03}:
$$
\begin{array}{ll}
&\Delta(B_n)=\{\epsi_i-\epsi_{i+1},1\leq i\leq n-1\}\cup\{\epsi_n\},\\
&\Delta(D_n)=\{\epsi_i-\epsi_{i+1},1\leq i\leq n-1\}\cup\{\epsi_{n-1}+\epsi_n\}.\\
&\Delta(C_n)=\{\epsi_i-\epsi_{i+1},1\leq i\leq n-1\}\cup\{2\epsi_n\}.
\end{array}
$$
The set of positive roots $\Phi^+\supset\Delta$ looks as follows:
$$
\begin{array}{ll}
&B_n^+ = \{\epsi_i\pm\epsi_j, 1\leq i < j\leq n\}\cup\{\epsi_i,
1\leq i\leq n\},\\
&D_n^+ = \{\epsi_i\pm\epsi_j, 1\leq i < j\leq n\},\\
&C_n^+ = \{\epsi_i\pm\epsi_j, 1\leq i < j\leq n\}\cup\{2\epsi_i,
1\leq i\leq n\}.
\end{array}
$$
Let
$$
m=\left\{\begin{array}{ll}2n+1,&\mbox{if }\Phi=B_n,\\
2n,&\mbox{if }\Phi=D_n\text{ or }C_n.
\end{array}\right.
$$
Denote by $\ut=\ut(\Phi)$ the subalgebra of $\mathfrak{gl}_m(\Fp_q)$
spanned by the vectors $e_{\alpha}$, $\alpha\in\Phi$, where
$$
\begin{array}{ll}
&e_{\epsi_i}=e_{0,i}-e_{-i,0},\quad 1\leq i\leq n,\\
&e_{\epsi_i-\epsi_j}=e_{j,i}-e_{-i,-j},\quad 1\leq i<j\leq n,\\
&e_{\epsi_i+\epsi_j}=\begin{cases}
e_{-j,i}+e_{-i,j},\text{ if }\Phi=C_n,\\
e_{-j,i}-e_{-i,j}\text{ otherwise},
\end{cases}
\quad 1\leq i<j\leq n.
\end{array}
$$
Here we numerate the rows and the columns of $m\times m$ matrix by the indices $$1,2,\ldots,n,0,-n,\ldots,-2,-1$$ (there is no index $0$ in the cases $D_n$ and $C_n$), and we denote by $e_{a, b}$ the usual elementary matrix. Of course, it is a maximal nilpotent subalgebra in the corresponding classical algebra $\gt= \gt(\Phi)$. In particular, $\dim\ut = |\Phi^+|$. \exam{Here we schematically drew the algebras $B_3$, $D_4$ and $C_4$. The convention is as follows: given $1\leq j<i$, the square $(i,j)$ (respectively, $(-i,j)$, $(0,j)$ and $(-j,j)$) corresponds to the root $\epsi_j-\epsi_i$ (respectively, $\epsi_j+\epsi_i$, $\epsi_j$ and $2\epsi_j$).
\begin{center}\small
$\mymatrix{
\Note{1}\lNote{1}\pho& \Note{2}\pho& \Note{3}\pho& \Note{0}\pho& \Note{-3}\pho& \Note{-2}\pho & \Note{-1}\pho\\
\lNote{2}\Top{2pt}\Rt{2pt} \pho & \pho& \pho& \pho& \pho& \pho & \pho\\
\lNote{3}\pho & \Top{2pt}\Rt{2pt} \pho& \pho& \pho& \pho& \pho & \pho\\
\lNote{0}\pho & \pho& \Top{2pt}\Rt{2pt} \pho& \pho& \pho& \pho & \pho\\
\lNote{-3}\pho & \pho& \Top{2pt}\Lft{2pt} & \Top{2pt}\Rt{2pt} \pho& \pho& \pho & \pho\\
\lNote{-2}\pho& \Top{2pt}\Lft{2pt} & \pho& \pho& \Top{2pt}\Rt{2pt} \pho& \pho & \pho\\
\lNote{-1}\Top{2pt} & \pho& \pho& \pho& \pho& \Top{2pt}\Rt{2pt}\pho & \pho\\
}\quad\mymatrix{
\lNote{1}\Note{1}\pho& \Note{2}\pho& \Note{3}\pho& \Note{4}\pho& \Note{-4}\pho& \Note{-3}\pho & \Note{-2}\pho & \Note{-1}\pho\\
\lNote{2}\Top{2pt}\Rt{2pt} \pho & \pho& \pho& \pho& \pho& \pho & \pho & \pho\\
\lNote{3}\pho & \Top{2pt}\Rt{2pt} \pho& \pho& \pho& \pho& \pho & \pho & \pho\\
\lNote{4}\pho & \pho& \Top{2pt}\Rt{2pt} \pho& \pho& \pho& \pho & \pho & \pho\\
\lNote{-4}\pho &\pho & \pho& \Top{2pt}\Lft{2pt}\Rt{2pt} & \pho& \pho& \pho & \pho\\
\lNote{-3}\pho &\pho& \Top{2pt}\Lft{2pt} & \pho& \Top{2pt}\Rt{2pt}\pho& \pho& \pho & \pho\\
\lNote{-2}\pho &\Top{2pt}\Lft{2pt} & \pho& \pho& \pho& \Top{2pt}\Rt{2pt} \pho& \pho & \pho\\
\lNote{-1}\Top{2pt} & \pho& \pho& \pho& \pho& \pho & \Top{2pt}\Rt{2pt}\pho &\pho\\
}\quad\mymatrix{
\lNote{1}\Note{1}\pho& \Note{2}\pho& \Note{3}\pho& \Note{4}\pho& \Note{-4}\pho& \Note{-3}\pho & \Note{-2}\pho & \Note{-1}\pho\\
\lNote{2}\Top{2pt}\Rt{2pt} \pho & \pho& \pho& \pho& \pho& \pho & \pho & \pho\\
\lNote{3}\pho & \Top{2pt}\Rt{2pt} \pho& \pho& \pho& \pho& \pho & \pho & \pho\\
\lNote{4}\pho & \pho& \Top{2pt}\Rt{2pt} \pho& \pho& \pho& \pho & \pho & \pho\\
\lNote{-4}\pho &\pho & \pho& \Top{2pt}\Bot{2pt}\Rt{2pt} & \pho& \pho& \pho & \pho\\
\lNote{-3}\pho &\pho& \Bot{2pt}\Rt{2pt} & \pho& \Top{2pt}\Rt{2pt}\pho& \pho& \pho & \pho\\
\lNote{-2}\pho &\Bot{2pt}\Rt{2pt} & \pho& \pho& \pho& \Top{2pt}\Rt{2pt} \pho& \pho & \pho\\
\lNote{-1}\Bot{2pt}\Rt{2pt}& \pho& \pho& \pho& \pho& \pho & \Top{2pt}\Rt{2pt}\pho &\pho\\
}
$
\end{center}
}

Furthermore, we define the functions 
$$
\begin{array}{ll}
&\col\colon\Phi^+\to\{1,\ldots,n\}\colon\col(\epsi_i\pm\epsi_j)=\col(\epsi_i)=\col(2\epsi_i)=i,\\
&\row\colon\Phi^+\to\{-n,\ldots,n\}\colon\row(\epsi_i\pm\epsi_j)=\mp j,\row(\epsi_i)=0,\row(2\epsi_i)=i.\\
\end{array}
$$
For arbitrary $-n+1\leq i\leq n-1$ and $1\leq j\leq n$, the sets
$$
\begin{array}{ll}
&R_i = R_i(\Phi) = \{\alpha\in\Phi^+\mid \row(\alpha)=i\},\\
&C_j = C_j(\Phi) = \{\alpha\in\Phi^+\mid \col(\alpha)=j\}\\
\end{array}
$$
are called the $i$th\emph{ row} and the $j$th\emph{ column}
$\Phi^+$ respectively. We introduce the mirror order on the 
set of indices $$1\prec2\prec\ldots\prec
n\prec0\prec-n\prec\ldots\prec-2\prec-1,$$ and the following total orders on $\Phi^+$:
\begin{equation*}
\alpha\prec\beta~\stackrel{\mathrm{def}}{\iff}~
\col(\beta)\prec\col(\alpha)~\mbox{or}~
\col(\beta)=\col(\alpha),~\row(\beta)\succ\row(\alpha).
\end{equation*}
For example, for $\Phi=B_6$ we have
$\epsi_2+\epsi_5\succ\epsi_2\succ\epsi_2-\epsi_4\succ\epsi_3-\epsi_6$.

For given $\Phi$ denote the maximal unipotent subgroup $\exp(\ut(\Phi))$ in the corresponding classical finite group $G$ by~$U(\Phi)$ (we will write just $\ut$ and $U$ in cases when $\Phi$ is fixed). In the sequel, we will assume everywhere that $p=\chara{\Fp_q}\gee n$. Under this assumption, the map
$$
\exp\colon\ut\to U\colon x\mapsto \sum_{i=0}^{m-1}\frac{x^i}{i!}
$$
is well-defined and is in fact a bijection (and also an isomorphism of algebraic varieties over $\overline{\Fp_q}$);
we denote the inverse map by $\ln$. Furthermore, the 
Backer--Campbell--Hausdorff formula claims that, for a Lie subalgebra $\at\subset\ut$ and arbitrary $u, v\in\at$, one has
\begin{equation*}
\exp(u)\exp(v)=\exp(u + v + \tau(u,
v)),
\end{equation*} where
$\tau(u, v)\in[\at, \at]$ (here $[\at, \at]=\langle[x, y], x,
y\in\at\rangle_{\Fp_q}$).

The group $U$ acts on its Lie algebra $\ut$ by the adjoint action; the dual action of $U$ on the $\Fp_q$-dual space $\ut^*$ is called \emph{coadjoint}. Using the non-degenerate form $\langle A, B\rangle=\frac{1}{2}\mathrm{tr}(AB)$ on $\mathfrak{gl}_m(\Fp_q)$, one can identify the dual space $\ut^*$ with the space $\ut^t$ (in this case, $e_{\alpha}^*=e_{\alpha}^t$ for any
$\alpha=2\epsi_i\in\Phi^+$ and $e^*_{\alpha}=e_{\alpha}^T/2$ otherwise). Under this identification, the coadjoint action has the following form:
$$
g.x = \mathrm{pr}(gxg^{-1}),~g\in U,~x\in\ut^*
$$
(here we denote the projection
$\mathfrak{gl}_m(\Fp_q)\to\ut^*$ along $\ut$ by $\mathrm{pr}$).
\defi{Let $D=\{\beta_1,\ldots,\beta_t\}\subset\Phi^+$ be an arbitrary subset of the positive roots, and
$\xi=(\xi_{\beta})_{\beta\in D}$ be a set of non-zero constants from $\Fp_q$. Put
$$
f=f_{D, \xi} = \sum_{\beta\in D}\xi_{\beta}e_{\beta}^*\in\ut^*. $$ We say that the coadjoint orbit $\Omega=\Omega_{D,
\xi}\subset\ut^*$ of the linear form $f$ is \emph{associated} with $D$, and $f$ is called the
\emph{canonical form} on $\Omega$.}

\defi{Let $\beta\in\Phi^+$. Roots $\alpha$, $\gamma\in\Phi^+$ are called
$\beta$-\emph{singular} if $\alpha+\gamma=\beta$. The set of all $\beta$-singular roots is denoted by $S(\beta)$ (see
\cite{Andre95ii} for $A_n$ and \cite{AndreNeto06} for
$B_n$, $D_n$).} It is easy to see that the singular roots look as follows:
\begin{equation*}
\begin{split}
&S(\epsi_i-\epsi_j)=\bigcup_{l=i+1}^{j-1}\{\epsi_i-\epsi_l,\epsi_l-\epsi_j\},~1\leq i < j\leq n,\\
&S(\epsi_i) = \bigcup_{l=i+1}^{n}\{\epsi_i-\epsi_l,\epsi_l\},~1\leq i\leq n,\\
&S(2\epsi_i) = \bigcup_{l=i+1}^{n}\{\epsi_i-\epsi_l,\epsi_i+\epsi_l\},~1\leq i\leq n,\\
&S(\epsi_i+\epsi_j)=\bigcup_{l=i+1}^{j-1}\{\epsi_i-\epsi_l,\epsi_l+\epsi_j\}
\cup\bigcup_{l=j+1}^n\{\epsi_i-\epsi_l,\epsi_j+\epsi_l\}\cup\\
&\bigcup_{l=j+1}^n\{\epsi_j+\epsi_l,\epsi_j-\epsi_l\}\cup S_{ij}
,~1\leq i<j\leq
n,\mbox{ where}\\
&S_{ij}=
\begin{cases}\{\epsi_i, \epsi_j\},&\mbox{if }\Phi=B_n,\\
\{\epsi_i-\epsi_j, 2\epsi_j\},&\mbox{if }\Phi=C_n\\
\varnothing,&\mbox{if }\Phi=D_n.
\end{cases}.
\end{split}\label{formula_sing_roots}
\end{equation*}

\nota{Note that, given $\beta\in\Phi^+$, 
the column of $\beta$ contains exactly one of the root from each pair of $\beta$-singular roots whose sum is $\beta$, except for the case $\beta=2\epsi_i$. Precisely, put
\begin{equation*}
S^+(\beta)=
\begin{cases}
\{\epsi_i+\epsi_j,~i<k\leqslant n\},&\text{ if }\beta=2\epsi_i,\\
S(\beta)\cap C_{\text{col}(\beta)}&\text{ otherwise},
\end{cases}
\end{equation*}
so that $S(\beta)=S^+(\beta)\cup S^-(\beta)$\label{nota_one_row_col}, where $S^-(\beta)=S(\beta)\backslash S^+(\beta)$.}

Next, we need to recall some necessary notations from \cite{IgnatevPetukhov25} which will be used in proofs of the main results of the paper.

Let $X$ be a subset of $\Phi^+$. Set $\ut_X:=\bigoplus_{\alpha\in X}\Fp e_\alpha$ with the (Lie) bracket given by the formula $$[e_\alpha, e_\beta]=\begin{cases}[e_\alpha, e_\beta],& {\rm~if~}\alpha+\beta\in X,\\0,&\text{if }\alpha+\beta\notin X\end{cases}$$ for all $\alpha, \beta\in X$
(note that if $\alpha, \beta, \alpha+\beta\in \Phi$, then $[e_\alpha, e_\beta]$ is proportional to $e_{\alpha+\beta})$. Such a bilinear bracket $[\cdot, \cdot ]$ do not always define a Lie algebra but it does define a Lie algebra under the assumption that $X$ is a $C$-pattern (a shorthand for ``combinatorial counterpart of a pattern subgroup'') or $C$-quattern (a shorthand for ``combinatorial counterpart of a quattern subgroup''), see the next definitions. 

\defi{We will say that $X$ is a {\it C-pattern}  if $X$ satisfies the following condition: $$\text{if }\alpha, \beta\in X\text{ and }\alpha+\beta\in\Phi^+\text{ then }\alpha+\beta\in X.$$
This condition is clearly equivalent to the statement that $\ut_X$ is a sub-Lie algebra of $\ut$.}
\defi{ We will say that $X$ is a {\it C-quattern} if $X=X_+\setminus X_-$ for two C-patterns $X_+, X_-$ such that $X_-\subset X_+$ and
$$\text{if }\alpha_-\in X_-,\alpha_+\in X_+~\text{ and }\alpha_-+\alpha_+\in X_+ \text{ then }\alpha_-+\alpha_+\in X_-.$$
This is clearly equivalent to the statement that $\ut_{X_-}$ is an ideal of the Lie algebra $\ut_{X_+}$ and in this case we have $\ut_X\cong \ut_{X_+}/\ut_{X_-}$ where the right hand side is a quotient of Lie algebra by its ideal and thus the left hand side is a Lie algebra.}
From now on we assume that $X$ is always a C-quattern. We denote by $U_X$ the unipotent group with Lie algebra $\ut_X$. Every C-pattern $X$ is a C-quattern for $X_+=X, X_-={\varnothing}$ and for every C-quattern $X$ the data $(\ut_X, [\cdot, \cdot])$ provides a structure of Lie algebra on $\ut_X$. 

Set 
$$Z(X):=\{\alpha\in X\mid (\alpha+X)\cap X={\varnothing}\}.$$ 
Then it is easy to verify that $\ut_{Z(X)}$ coincides with the center of $\ut_X$. In particular, this implies that $Z(X)\ne\varnothing$ for any C-quattern $X$.

\defi{Pick $f\in \ut_X^*$ and consider $Z\subset Z(X)$. We will say that $f$ is {\it $Z$-saturated} if $f(e_\alpha)\ne0$ for all $\alpha\in Z$. We will say that $f$ is {\it saturated} if $f$ is $Z(X)$-saturated. 
We denote the variety of $Z$-saturated linear forms by $\ut^*_{X; Z}$.}
It is clear that each $\ut_{X; Z}^*$ is $U_X$-stable and thus is a collection of coadjoint orbits. We map $\ut^*_{X; Z}$ to $\ut^*$ by the formula $e_\alpha^*\mapsto e_\alpha^*$ and denote by $\underline{\ut}^*_{X; Z}$ the image of this map. 
We will frequently identify $\ut_{X; Z}^*$ with $\underline{\ut}_{X; Z}^*$. 

\defi{Let $Z\subset  Z(X)$ be a subset and let $Y\subset \ut_{X; Z}^*$ be a subvariety. 
We say that $Y$ is a {\em set-section} for $\ut_{X; Z}^*$ if $Y$ intersects each $U_X$-orbit of $\ut_{X; Z}^*$ in exactly one point.}
It is clear that a set-section is a section of the native quotient map $\ut_{X; Z}^*\to \ut_{X; Z}^*/U_X$  but only in set-theoretic terms. 
Also note that all the set-sections which we will consider explicitly will be unions of the form $V(S_1)\sqcup V(S_2)\sqcup\ldots$ for some subsets $S_1, S_2, \ldots\subset\Phi^+$.

\sect{Subregular orbits of group of type $C_n$}\label{sect:sea_battle_pol}

Before formulating the main theorem of this paragraph, we need to define some subsets of the set positive roots. Put $$D_\text{reg}=\{2\epsi_1,2\epsi_2,\ldots,2\epsi_n,\epsi_1+\epsi_2,\epsi_2+\epsi_3,\ldots,\epsi_{n-1}+\epsi_n\}\subset\Phi^+,$$
and $$D'=\{\epsi_1+\epsi_3,\epsi_2+\epsi_4,\ldots,\epsi_{n-2}+\epsi_n,\epsi_{n-1}-\epsi_{n}\}.$$
Now, we want to introduce the definition of a subregular subset of roots. First, we pick a root $\alpha$ for from the set $D'$. Next, we choose exact one root $\beta$ from the set $D''=\{2\epsi_{i+1},\gamma\}$, where $\gamma$ is the closest root to $\alpha$, which is smaller than $\alpha$, from the set $D'$ with respect to the order $\prec$ and $i=\col(\alpha)$ (in the case, when $\alpha=\epsi_{n-1}-\epsi_n$, the set $D''=\{2\epsi_n\}$). At the next step, we add a set of roots $\mathcal{D}$ from the set $D_\text{reg}$ to the set $D''$ such that $\mathcal{D}\cup\{\alpha,\beta\}$ is maximal orthogonal set of roots (we call a set of roots orthogonal if any two roots from this set are orthogonal) under the assumption that for $\alpha=\epsi_{n-2}+\epsi_n$ and $\beta=\epsi_{n-1}-\epsi_n$ or $\alpha=\epsi_{n-1}$ and $\beta=2\epsi_n$ we consider $\alpha$ and $\beta$ orthogonal. Finally, we add the root $\epsi_{i}-\epsi_{i+1}$ to the set $\mathcal{D}\cup\{\alpha,\beta\}$ for $\alpha=\epsi_{i-1}+\epsi_{i+1}$. A subset of roots constructed as above we will call subregular and denote $D_\text{sreg}$.  

As in Section \ref{main_definitions}, we will consider a set $\xi$ of constants from $\Fp_q^{\times}$; however in this section, it will be defined with a minor correction. More precisely, put $\xi=(\xi_{\beta})_{\beta\in D_\text{sreg}}$ for any $D_\text{sreg}$, but allow $\xi_{2\epsi_n}$ and $\xi_{\epsi_i-\epsi_{i+1}}$ to be zero for $i$ from 1 to $n-1$. Now we are ready to formulate and prove the following result.
\exam{Here we schematically drew some examples of elements $f_{D,\xi}$ for some $D=D_\text{sreg}$ for $\Phi=C_4$. The symbol $\otimes$ stands for nonzero constants and the symbol $\square$ for an arbitrary constants $\xi_{\beta}$ for $\beta\in D_\text{sreg}$.
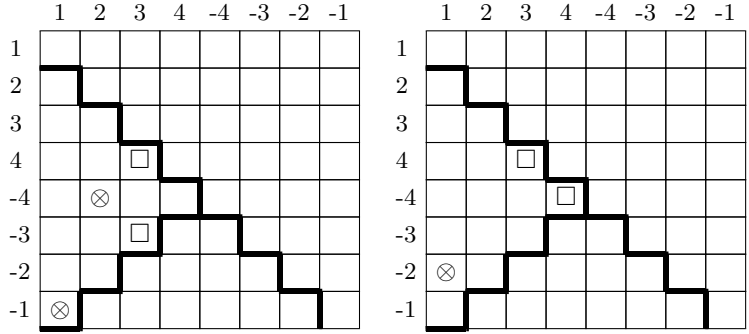
\begin{center}\small
 $\mymatrix{
\lNote{1}\Note{1}\pho& \Note{2}\pho& \Note{3}\pho& \Note{4}\pho& \Note{-4}\pho& \Note{-3}\pho & \Note{-2}\pho & \Note{-1}\pho\\
\lNote{2}\Top{2pt}\Rt{2pt} \pho & \pho& \pho& \pho& \pho& \pho & \pho & \pho\\
\lNote{3}\pho & \Top{2pt}\Rt{2pt} \pho& \pho& \pho& \pho& \pho & \pho & \pho\\
\lNote{4}\pho & \pho& \Top{2pt}\Rt{2pt}\square& \pho& \pho& \pho & \pho & \pho\\
\lNote{-4}\pho &\otimes & \pho& \Top{2pt}\Bot{2pt}\Rt{2pt} & \pho& \pho& \pho & \pho\\
\lNote{-3}\pho &\pho& \Bot{2pt}\Rt{2pt}\square& \pho& \Top{2pt}\Rt{2pt}\pho& \pho& \pho & \pho\\
\lNote{-2}\pho &\Bot{2pt}\Rt{2pt} & \pho& \pho& \pho& \Top{2pt}\Rt{2pt} \pho& \pho & \pho\\
\lNote{-1}\Bot{2pt}\Rt{2pt}\otimes& \pho& \pho& \pho& \pho& \pho & \Top{2pt}\Rt{2pt}\pho &\pho\\
}
\quad\mymatrix{
\lNote{1}\Note{1}\pho& \Note{2}\pho& \Note{3}\pho& \Note{4}\pho& \Note{-4}\pho& \Note{-3}\pho & \Note{-2}\pho & \Note{-1}\pho\\
\lNote{2}\Top{2pt}\Rt{2pt} \pho & \pho& \pho& \pho& \pho& \pho & \pho & \pho\\
\lNote{3}\pho & \Top{2pt}\Rt{2pt} \pho& \pho& \pho& \pho& \pho & \pho & \pho\\
\lNote{4}\pho & \pho& \Top{2pt}\Rt{2pt}\square& \pho& \pho& \pho & \pho & \pho\\
\lNote{-4}\pho &\pho & \pho& \Top{2pt}\Bot{2pt}\Rt{2pt}\square& \pho& \pho& \pho & \pho\\
\lNote{-3}\pho &\pho& \Bot{2pt}\Rt{2pt} & \pho& \Top{2pt}\Rt{2pt}\pho& \pho& \pho & \pho\\
\lNote{-2}\otimes &\Bot{2pt}\Rt{2pt} & \pho& \pho& \pho& \Top{2pt}\Rt{2pt} \pho& \pho & \pho\\
\lNote{-1}\Bot{2pt}\Rt{2pt}& \pho& \pho& \pho& \pho& \pho & \Top{2pt}\Rt{2pt}\pho &\pho\\
}$
\end{center}
At the left picture we have $D=\{2\epsi_1,\epsi_2+\epsi_4,2\epsi_3,\epsi_3+\epsi_4\}$ and $D=\{\epsi_1+\epsi_2,\epsi_3+\epsi_4,2\epsi_4\}$ at the right picture.
}
\theop{Let\label{theo_orb_max_dim} $D=D_\text{sreg}\subset\Phi^+$ be a subregular subset of the positive roots and $\xi$ be as above\textup{,} then put 
\begin{equation*}
f=f_{D, \xi} = \sum_{\beta\in D}\xi_{\beta}e_{\beta}^*\in\ut^*.
\end{equation*}
Then the orbit of $f$ has pre-maximal dimension. Conversely\textup{,} each orbit of pre-maximal dimension contains exact one such a linear form.}{We will proceed by induction on $n$. The base of induction for $C_2$ and $C_3$ follows directly from \cite{GoodwinMoschRoehrle16}. Firstly we want to proof that, for given $f=(\xi_{\beta})_{\beta\in\Phi^+}\in\ut^*$ and $\delta\prec\epsi_1+\epsi_3$, where $\delta$ is the biggest root with nonzero constant $\xi_{\delta}$ with respect to the order $\prec$, the orbit of $f$ has a dimension which is lower than pre-maximal. Let us assume that $\delta=\epsi_1+\epsi_4$, i.e., $\xi_{2\epsi_1}=\xi_{\epsi_1+\epsi_2}=\xi_{\epsi_1+\epsi_3}=0$. Consider the quattern $Y=\Phi^+\backslash\{2\epsi_1,\epsi_1+\epsi_2,\epsi_1+\epsi_3\}$ and denote by $\wt f$ the projection of $f$ to $\ut_{Y;Z}$ for $Z=\langle e_{\delta}\rangle\subset Z(Y)$. The dimension of the orbit of the form $f$ can be calculated as $\text{rank}(\beta_{f})$ (recall that $\beta_f$ is a quadratic form on $\ut$ defined by $\beta_f(x,y)=f([x,y])$). We want to compare the ranks of the quadratic forms $\beta_f$ and $\beta_{\wt f}$ defined on $\ut$ and $\ut_{Y;Z}$ respectively. Let us denote the matrix of the form $\beta_f$ by $A$ and the matrix of the form $\beta_{\wt f}$ by $\wt A$ and write the ideal spanned by the vectors $e_{2\epsi_1},e_{\epsi_1+\epsi_2},e_{\epsi_1+\epsi_3}$ as $I$ (columns and rows are numbered by the roots from $\Phi^+$ which are sorted in the order $\prec$). Hence, the matrix $A$ consist of the numbers $\beta_f(e_{\alpha},e_{\beta})\in\Fp_q$ on $\alpha$-th, $\beta$-th places. If $[e_{\alpha},e_{\beta}]\notin I$ for some $\alpha,\beta\in\Phi^+$, then $\beta_f(e_{\alpha},e_{\beta})=\beta_{\wt f}(e'_{\alpha},e'_{\beta})$ (here $e'_{\alpha}$ is the class of the vector $e_{\alpha}$), if $[e_{\alpha},e_{\beta}]\in I$, then $\beta_f(e_{\alpha},e_{\beta})=\beta_{\wt f}(e'_{\alpha},e'_{\beta})=0$ because $f|_{I}=0$ and $[e'_{\alpha},e'_{\beta}]=0$ and finally, if $e_{\alpha}$ of $e_{\beta}$ belong to $I$, we have $\beta_f(e_{\alpha},e_{\beta})=0$ and the corresponding row and column in $\wt A$ does not even exist. Thus, the matrix $\wt A$ is obtained from $A$ by adding some strings and columns of zeroes, so their ranks are equal. 


Hence, we can consider the form $\wt f$ and its orbit in $\ut^*_{Y;Z}$ instead of $f$ and its orbit in $\ut^*$.
Let $Y'$ be defined as the set $Y\backslash\{\delta_1,\cdots,\delta_{2n-5},\beta_1\cdots,\beta_{2n-5}\}$, where $\delta_i\in S^+(\delta)$ and $\beta_i\in S^-(\delta)$. Consider the map $\vfi:\ut^*_{Y}\to\ut^*=\ut^*_X$ defined as follows for any quattern $Y$: for $e_{\tau}^*\in\ut_Y^*$,
$$\vfi(e^*_{\tau})=e^*_{\tau}\in\ut^*.$$ According to \cite[Proposition 5.11]{IgnatevPetukhov25},
the set $Y'$ is a quattern and there is a form $f'\in\ut^*_{Y';Z}$ such that $\vfi(f')$ lying on the orbit of the form $\wt f$ in $\ut^*_{Y;Z}$ and $$\dim\Omega_{\wt f}=\dim\Omega_{f'}+|\{\delta_1,\ldots,\delta_{2n-5},\beta_1,\ldots,\beta_{2n-5}\}|=\dim\Omega_{f'}+2(2n-5).$$ Now, we want to examine the maximal possible dimension of orbit in $\ut^*_{Y';Z}$. Let us assume that $f'\in\ut^*_{Y';Z'}$, where $Z'=Z\cup\langle e_{2\epsi_2}\rangle$ (i.e. $f'(e_{2\epsi_2})\neq0$). Put $Y''=Y'\backslash\{\delta'_1,\ldots,\delta'_{n-3},\beta'_1,\ldots,\beta'_{n-3}\}$ for $\delta'_i\in S^+(2\epsi_2)$ and $\beta'_i\in S^-(2\epsi_2)$ but in sense of the quattern $Y'$. According to the Proposition 5.11, the set $Y''$ is a quattern, and there exist $f''\in\ut^*_{Y'';Z'}$ such that $\vfi'(f'')=f'$, where $\vfi'=\vfi|_{Y''}$ such that  
$$\dim\Omega_{f'}=\dim\Omega_{f''}+|\{\delta'_1,\ldots,\delta'_{n-3},\beta'_1,\ldots,\beta'_{n-3}\}|=\dim\Omega_{f''}+2(n-3).$$

It is easy to check that the quotient algebra $\ut^*_{Y''}$ is isomorphic to the sum of its ideals $\ut(C_{n-3})\oplus\langle e_{2\epsi_2}\rangle\oplus\langle e_{\epsi_2-\epsi_4}\rangle\oplus\langle e_{\epsi_3-\epsi_4}\rangle$. It follows that $\dim\Omega_{f''}$ in $\ut^*_{Y'';Z'}$ is no more than $(n-3)(n-4)$, the maximal possible dimension of orbit in $\ut^*(C_{n-3})$. So, $\dim\Omega_{f'}$ in $\ut^*_{Y';Z'}$ is no more than $(n-3)(n-4)+2(n-3)$. It is left to note that $\ut^*_{Y';Z'}$ is an open subspace of $\ut^*_{Y';Z}$ as well as the set of orbits of maximal dimension, so its intersection with the set of orbits of maximal dimension is non-empty, so the maximal dimension of an orbit in $\ut^*_{Y',Z}$ is also no more than $(n-3)(n-4)+2(n-3)$. 
Hence, we obtain that $$\dim\Omega_{\wt f}\leqslant (n-3)(n-4)+2(n-3)+2(2n-5)=n(n-1)-4,$$
what is less than pre-maximal dimension of the orbit in $\ut^*$.

It follows that the maximal possible dimension in $\ut^*_{Y,Z}$ is less than pre-maximal. Using the previous arguments, we obtain that the maximal possible dimension in $\ut^*_{Y}$ is also less than pre-maximal. So, if the root $\delta$ is bigger than $\epsi_1+\epsi_3$, we can consider the orbit of $f$ as orbit in the quotient algebra $\ut^*_{Y}$. It concludes the poof of the first part.



Now, we know that a linear form $f=(\xi_{\beta})_{\beta\in\Phi^+}$ having the orbit of maximal dimension has to have nonzero constant $\xi_\delta$, where $\delta$ is the biggest root such that $\xi_{\delta}\ne 0$ with $\delta=2\epsi_1$, $\epsi_1+\epsi_2$ or $\epsi_1+\epsi_3$. We aim to present a set-section of the set of orbits of pre-maximal dimension. Let us consider the cases when $\delta=2\epsi_1$ and $\delta=\epsi_1+\epsi_3$, the case $\delta=\epsi_1+\epsi_2$ is similar to the case, when $\delta=2\epsi_1$. We proceed by induction on $n$. Firstly, put $X=\Phi^+$, $\delta=2\epsi_1$, $Z=\langle\delta\rangle= Z(X)$, and $$Y=X\backslash\{\delta_1,\dots,\delta_{n-1},\beta_1,\dots,\beta_{n-1}\}$$ for $\delta_i=\epsi_1-\epsi_i$ and $\beta_i=\epsi_1+\epsi_i$. For an arbitrary form $f$ as above, we have $f\in\ut_{X;Z}$ by the definition of $\ut_{X;Z}$. Note that the set $Y$ is a quattern and it is isomorphic to $\Phi(C_{n-1})\cup\delta$, so $\ut^*_{Y;Z}$ is the direct sum of $\ut^*(C_{n-1})\oplus\langle\delta\rangle$ as ideals. The Proposition 5.11 gives us that for a set-section $\wt K$ of $\ut^*_{Y;Z}$, the set $K=\vfi(\wt K)$ is a set-section for $\ut^*_{X;Z}$, moreover if $\wt f\in\wt K$ for $f=\vfi(\wt f)\in K$, then 
$$\dim\Omega_{f}=\dim\Omega_{\wt f}+|\{\delta_1,\ldots,\delta_{n-1},\beta_1,\ldots,\beta_{n-1}\}|=\dim\Omega_{\wt f}+2(n-1).$$
It follows that dimension of $f$ is pre-maximal (i.e. $n(n-1)-2$) in $\ut^*_{X;Z}$ if and only if dimension of $\wt f$ is pre-maximal (i.e $(n-1)(n-2)-2$) in $\ut^*_{Y;Z}$. It means that we can apply the induction assumption to a set-section $\wt K$ of $\ut^*_{Y;Z}$. It concludes the proof in that case.

Now, put $X=\Phi^+\backslash\{2\epsi_1,\epsi_1+\epsi_2\}$, $\delta=\epsi_1+\epsi_2$, $Z=\langle\delta\rangle\subset Z(X)$, and $$Y=X\backslash\{\delta_1,\dots,\delta_{2n-4},\beta_1,\dots,\beta_{2n-4}\}$$ for $\delta_i$ from $S^+(\delta)$ and $\beta_i$ from $S^-(\delta)$. For an arbitrary form $f$ as above, we have $f\in\ut_{X;Z}$ by the definition of $\ut_{X;Z}$. Note that the set $Y$ is a quattern and it is isomorphic to $\Phi(C_{n-2})\cup\delta\cup\epsi_2-\epsi_3$, so $\ut^*_{Y;Z}$ is the direct sum of $\ut^*(C_{n-2})\oplus\langle\delta\rangle\oplus\langle\epsi_2-\epsi_3\rangle$ as ideals. The Proposition 5.11 gives us that for a set-section $\wt K$ of $\ut^*_{Y;Z}$, the set $K=\vfi(\wt K)$ is a set-section for $\ut^*_{X;Z}$, moreover if $\wt f\in\wt K$ for $f=\vfi(\wt f)\in K$, then 
$$\dim\Omega_{f}=\dim\Omega_{\wt f}+|\{\delta_1,\ldots,\delta_{n-1},\beta_1,\ldots,\beta_{n-1}\}|=\dim\Omega_{\wt f}+2(2n-4).$$
It follows that dimension of $f$ is pre-maximal (i.e. $n(n-1)-2$) in $\ut^*_{X;Z}$ if and only if dimension of $\wt f$ is maximal (i.e $(n-2)(n-3)$) in $\ut^*_{Y;Z}$. It means that we can apply a classification of the regular orbits in symplectic groups of type $C_n$ (see, for example, \cite{Venchakov}) to a set-section $\wt K$ of $\ut^*_{Y;Z}$. It concludes the proof in that case.
}
As a corollary, we can proof the Isaacs' conjecture in that case.
\corop{The number\label{coro:Isaacs_C} of subregular orbits of group of type $C_n$ over the finite field $\Fp_q$ is the polynomial from $v=q-1$ with integer non-negative coefficients.}{From the proof of the theorem above we see that such a number $S_n$ is given by the recurrent condition $$S_n=vS_{n-1}+vS_{n-2}+v(v+1)S'_{n-2},$$ where $S'_{n}$ is the expression for the number of regular orbits of group of type $C_n$ (see, for example \cite{Venchakov}). The base of induction gives us that $S_2=(v+1)^2$ and $S_3=2v(v+1)^2$. So, assuming that $S_{n-1}$, $S_{n-2}$ and $S'_{n-2}$ are polynomial from $v$ with integer non-negative coefficients, we obtain the required statement.}
\sect{Regular orbits of groups of type $B_n$ and $D_n$}

As in the previous section, we need to introduce some notations. More precisely, for $\Phi=B_n$ or $D_n$, put
$$D^{B_n}=\{\epsi_1+\epsi_2,\epsi_3+\epsi_4,\ldots,\epsi_{k-1}+\epsi_k,\epsi_1-\epsi_2,\epsi_3-\epsi_4,\ldots,\epsi_{k-1}-\epsi_k\}\subset\Phi^+,$$
where $k=n-2$ if $n$ is even and $k=n-1$ if $n$ is odd, and
$$D^{D_n}=\{\epsi_1+\epsi_2,\epsi_3+\epsi_4,\ldots,\epsi_{k-1}+\epsi_k,\epsi_1-\epsi_2,\epsi_3-\epsi_4,\ldots,\epsi_{k-1}-\epsi_k\}\subset\Phi^+,$$
where $k=n$ if $n$ is even and $k=n-1$ if $n$ is odd.

Next, in the case, when $\Phi=D_n$ we will denote the set $D^{D_n}$ by $D_\text{reg}$. And in the case $\Phi=B_n$, we will denote by $D_\text{reg}$ a one of the sets $D^{B_n}\cup\{\epsi_{n-1}+\epsi_n,\epsi_{n-1}-\epsi_n\}$ or $D^{B_n}\cup\{\epsi_{n-1}\}$ if $n$ is even, and $D^{B_n}\cup\{\epsi_n\}$ if $n$ is odd. Such a sets $D_\text{reg}$ of positive roots we will call regular. A set $\xi=(\xi_{\beta})_{\beta\in D_\text{reg}}$ is defined as above, except that $\xi_{\epsi_i-\epsi_{i+1}}$ for $i$ from 1 to $n-1$ and $\xi_{\epsi_n}$ are allowed to be zero.

\exam{Here we schematically drew some examples of elements $f_{D,\xi}$ for some $D=D_\text{reg}$ for $\Phi=B_4$ and $D_4$. The symbol $\otimes$ stands for nonzero constants and the symbol $\square$ for an arbitrary constants $\xi_{\beta}$ for $\beta\in D_\text{reg}$.
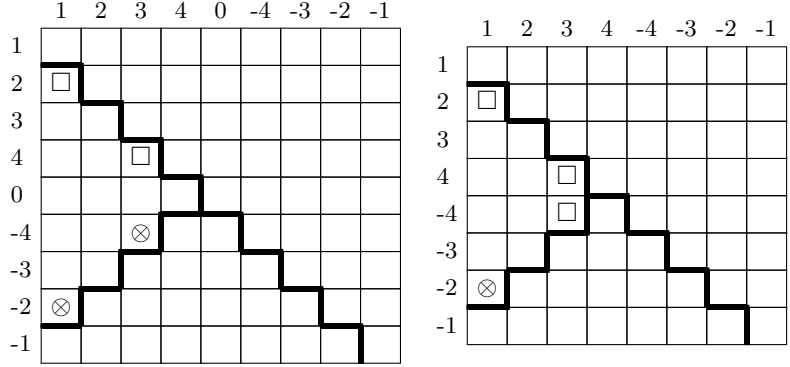
\begin{center}\small
$\mymatrix{
\lNote{1}\Note{1}\pho& \Note{2}\pho& \Note{3}\pho& \Note{4}\pho& \Note{0}\pho& \Note{-4}\pho& \Note{-3}\pho & \Note{-2}\pho & \Note{-1}\pho\\
\lNote{2}\Top{2pt}\Rt{2pt}\square & \pho& \pho& \pho& \pho& \pho& \pho & \pho & \pho\\
\lNote{3}\pho & \Top{2pt}\Rt{2pt} \pho& \pho& \pho& \pho& \pho& \pho & \pho & \pho\\
\lNote{4}\pho & \pho& \Top{2pt}\Rt{2pt}\square& \pho& \pho& \pho & \pho& \pho & \pho\\
\lNote{0}\pho & \pho& \pho& \Top{2pt\Rt{2pt}}\pho& \pho& \pho & \pho& \pho & \pho\\
\lNote{-4}\pho &\pho & \otimes& \Lft{2pt}\Top{2pt}& \Top{2pt}\Rt{2pt}& \pho& \pho& \pho & \pho\\
\lNote{-3}\pho &\pho& \Top{2pt}\Lft{2pt} & \pho& \pho& \Top{2pt}\Rt{2pt}& \pho& \pho & \pho\\
\lNote{-2}\otimes &\Top{2pt}\Lft{2pt} & \pho& \pho& \pho& \pho& \Top{2pt}\Rt{2pt}& \pho & \pho\\
\lNote{-1}\Top{2pt} & \pho& \pho& \pho& \pho& \pho & \pho& \Top{2pt}\Rt{2pt}& \pho\\
}\quad \mymatrix{
\lNote{1}\Note{1}\pho& \Note{2}\pho& \Note{3}\pho& \Note{4}\pho& \Note{-4}\pho& \Note{-3}\pho & \Note{-2}\pho & \Note{-1}\pho\\
\lNote{2}\Top{2pt}\Rt{2pt}\square& \pho& \pho& \pho& \pho& \pho & \pho & \pho\\
\lNote{3}\pho & \Top{2pt}\Rt{2pt} \pho& \pho& \pho& \pho& \pho & \pho & \pho\\
\lNote{4}\pho & \pho& \square\Top{2pt}\Rt{2pt}& \pho& \pho& \pho & \pho & \pho\\
\lNote{-4}\pho &\pho & \square\Bot{2pt}\Rt{2pt}& \Rt{2pt}\Top{2pt}& \pho& \pho& \pho & \pho\\
\lNote{-3}\pho & \Bot{2pt}\Rt{2pt}& \pho& \pho& \Top{2pt}\Rt{2pt}\pho& \pho& \pho & \pho\\
\otimes\lNote{-2}\Bot{2pt}\Rt{2pt} &\pho & \pho& \pho& \pho& \Top{2pt}\Rt{2pt} \pho& \pho & \pho\\
\lNote{-1}& \pho& \pho& \pho& \pho& \pho & \Top{2pt}\Rt{2pt}\pho &\pho\\
}
$
\end{center}
At the both pictures we drew the set $D=\{\epsi_1+\epsi_2,\epsi_1-\epsi_2,\epsi_3+\epsi_4,\epsi_3-\epsi_4\}$}

\theop{Let\label{theo_orb_max_dim_BD} $D=D_\text{reg}\subset\Phi^+$ be a regular subset of the positive roots for $\Phi=B_n$ or $D_n$ and $\xi$ be as above\textup{,} then put 
\begin{equation*}
f=f_{D, \xi} = \sum_{\beta\in D}\xi_{\beta}e_{\beta}^*\in\ut^*.
\end{equation*}
Then the orbit of $f$ has maximal dimension. Conversely\textup{,} each orbit of maximal dimension contains exact one such a linear form.}{As in the previous section, we will proceed by induction on $n$. The base of induction for $D_3$, $D_4$ and $B_2$, $B_3$ follows directly from \cite{GoodwinMoschRoehrle16}. The scheme of proof will be the same as in the previous theorem. Firstly, we will proof that for given $f=(\xi_{\beta})_{\beta\in\Phi^+}\in\ut^*$ and $\delta\prec\epsi_1+\epsi_2$, where $\delta$ is the biggest root with nonzero constant $\xi_{\delta}$ with respect to the order $\prec$, the orbit of $f$ has a dimension which is lower than maximal. Next, we will consider the case, when $\delta=\epsi_1+\epsi_2$, which is the only left case.

Let us move to the first statement. We assume that $\delta=\epsi_1+\epsi_3$. Arguing as in the previous theorem, we can conclude that we can consider the form $\wt f$ and it's orbit in $\ut^*_{Y;Z}$ instead of $f$ and it's orbit in $\ut^*$ for $Y=\Phi^+\backslash\{\epsi_1+\epsi_2\}$ and $Z=\langle e_{\epsi_1+\epsi_3}\rangle$. Let $Y'$ be defined to be the set $Y\backslash\{\delta_1,\cdots,\delta_k,\beta_1\cdots,\beta_k\}$, where $\delta_i\in S^+(\delta)\backslash\{\epsi_1-\epsi_3\}$ and $\beta_i\in S^-(\delta)\backslash\{\epsi_2+\epsi_3\}$ and $k=2n-5$ for $\Phi=B_n$ or $k=2n-6$ for $\Phi=D_n$. Consider the map $\vfi:\ut^*_{Y}\to\ut^*=\ut^*_X$ defined as follows for any quattern $Y$: for $e_{\tau}^*\in\ut_Y^*$,
$$\vfi(e^*_{\tau})=e^*_{\tau}\in\ut^*.$$ According to \cite[Proposition 5.11]{IgnatevPetukhov25},
the set $Y'$ is a quattern and there is a form $f'\in\ut^*_{Y';Z}$ such that $\vfi(f')$ lying on the orbit of the form $\wt f$ in $\ut^*_{Y;Z}$ and $$\dim\Omega_{\wt f}=\dim\Omega_{f'}+|\{\delta_1,\ldots,\delta_k,\beta_1,\ldots,\beta_k\}|=\dim\Omega_{f'}+2k.$$ Now, we want to examine the maximal possible dimension of orbit in $\ut^*_{Y';Z}$. Let us assume that $f'\in\ut^*_{Y';Z'}$, where $Z'=Z\cup\langle e_{\epsi_1-\epsi_3}\rangle$ (i.e. $f'(e_{\epsi_1-\epsi_3})\neq0$). Put $Y''=Y'\backslash\{\epsi_1-\epsi_2,\epsi_2-\epsi_3\}$. According to the Proposition 5.11, the set $Y''$ is a quattern, and there exist $f''\in\ut^*_{Y'';Z'}$ such that $\vfi'(f'')=f'$, where $\vfi'=\vfi|_{Y''}$ such that  
$$\dim\Omega_{f'}=\dim\Omega_{f''}+2.$$

It is easy to check that the quotient algebra $\ut^*_{Y''}$ is isomorphic to the sum of it's ideals $\ut(\Phi_{n-2})\oplus\langle e_{\epsi_2+\epsi_3}\rangle\oplus\langle e_{\epsi_1-\epsi_3}\rangle\oplus\langle e_{\epsi_1+\epsi_3}\rangle$. It follows that $\dim\Omega_{f''}$ in $\ut^*_{Y'';Z'}$ is no more than $l=(n-2)(n-3)$ for $\Phi=B_n$ and $l=(n-2)(n-4)$ or $(n-3)^2$ for $\Phi=D_n$ and $n$ is even or odd, respectively (what is the maximal possible dimension of orbit in $\ut^*(\Phi_{n-2})$). So, $\dim\Omega_{f'}$ in $\ut^*_{Y';Z'}$ is no more than $l+2$. It is left to note that $\ut^*_{Y';Z'}$ is an open subset of $\ut^*_{Y';Z}$ as well as the set of orbits of maximal dimension, so it's intersection with the set of orbits of maximal dimension is non-empty, so the maximal dimension of an orbit in $\ut^*_{Y',Z}$ is also no more than $l+2$. 
Hence, we obtain that $$\dim\Omega_{\wt f}\leqslant l+2+2k=\wt k-2,$$
where $\wt k$ is the maximal dimension of orbit in the corresponding root system.

It follows that the maximal possible dimension in $\ut^*_{Y,Z}$ is less than maximal. Using the previous arguments, we obtain that the maximal possible dimension in $\ut^*_{Y}$ is also less than maximal. So, if the root $\delta$ is bigger than $\epsi_1+\epsi_2$, we can consider the orbit of $f$ as orbit in the quotient algebra $\ut^*_{Y}$. It concludes the poof of the first part.

It is left to consider the case, when $\delta=\epsi_1+\epsi_2$ for each root system. Put $\Phi=D_n$, $X=\Phi^+$ and $Z=\langle e_{\epsi_1+\epsi_2}\rangle$. As previously, we consider the set of positive roots $$Y=X\backslash\{\delta_1,\dots,\delta_{2n-4},\beta_1,\dots,\beta_{2n-4}\}$$
for $\delta_i$ from $S^+(\delta)$ and $\beta_i$ from $S^-(\delta)$. For an arbitrary form $f$ as above, we have $f\in\ut_{X;Z}$ by the definition of $\ut_{X;Z}$. Note that the set $Y$ is a quattern and it is isomorphic to $\Phi^+(D_{n-2})\cup\{\delta\}\cup\{\epsi_1-\epsi_2\}$, so $\ut^*_{Y;Z}$ is the direct sum of $\ut^*(D_{n-2})\oplus\langle e_\delta\rangle\oplus\langle e_{\epsi_1-\epsi_2}\rangle$ as ideals. The Proposition 5.11 gives us that for a set-section $\wt K$ of $\ut^*_{Y;Z}$, the set $K=\vfi(\wt K)$ is a set-section for $\ut^*_{X;Z}$, moreover if $\wt f\in\wt K$ for $f=\vfi(\wt f)\in K$, then 
$$\dim\Omega_{f}=\dim\Omega_{\wt f}+|\{\delta_1,\ldots,\delta_{2n-4},\beta_1,\ldots,\beta_{2n-4}\}|=\dim\Omega_{\wt f}+2(2n-4).$$
Simple calculation shows us that dimension of $f$ is maximal (i.e. $n(n-2)$ or $(n-1)^2$ if $n$ is even or odd, respectively) in $\ut^*_{X;Z}$ if and only if dimension of $\wt f$ is maximal (i.e $(n-2)(n-4)$ or $(n-3)^2$ if $n$ is even or odd, respectively) in $\ut^*_{Y;Z}$. It means that we can apply the induction assumption to a set-section $\wt K$ of $\ut^*_{Y;Z}$. It concludes the proof in that case.

At last, put $\Phi=B_n$, $X=\Phi^+$ and $Z=\langle e_{\epsi_1+\epsi_2}\rangle$. As previously, we consider the set of positive roots $$Y=X\backslash\{\delta_1,\dots,\delta_{2n-3},\beta_1,\dots,\beta_{2n-3}\}$$
for $\delta_i$ from $S^+(\delta)$ and $\beta_i$ from $S^-(\delta)$. For an arbitrary form $f$ as above, we have $f\in\ut_{X;Z}$ by the definition of $\ut_{X;Z}$. Note that the set $Y$ is a quattern and it is isomorphic to $\Phi^+(B_{n-2})\cup\{\delta\}\cup\{\epsi_1-\epsi_2\}$, so $\ut^*_{Y;Z}$ is the direct sum of $\ut^*(B_{n-2})\oplus\langle e_\delta\rangle\oplus\langle e_{\epsi_1-\epsi_2}\rangle$ as ideals. The Proposition 5.11 gives us that for a set-section $\wt K$ of $\ut^*_{Y;Z}$, the set $K=\vfi(\wt K)$ is a set-section for $\ut^*_{X;Z}$, moreover if $\wt f\in\wt K$ for $f=\vfi(\wt f)\in K$, then 
$$\dim\Omega_{f}=\dim\Omega_{\wt f}+|\{\delta_1,\ldots,\delta_{2n-3},\beta_1,\ldots,\beta_{2n-3}\}|=\dim\Omega_{\wt f}+2(2n-3).$$
Simple calculation shows us that dimension of $f$ is maximal (i.e. $n(n-1)$) in $\ut^*_{X;Z}$ if and only if dimension of $\wt f$ is maximal (i.e $(n-2)(n-3)$) in $\ut^*_{Y;Z}$. It means that we can apply the induction assumption to a set-section $\wt K$ of $\ut^*_{Y;Z}$. It concludes the proof in that case.}

\nota{We want to notice that this method in this form can not be applied to obtain a classification of subregular orbits for groups of type $B_n$ and $D_n$. Consider the case of $D_4$, for example. In the situation when the root $\delta$ from the proof above is $\epsi_1+\epsi_4$ we can not pick any roots $\delta_i$ and $\beta_i$ to use the Proposition 5.11. We face the similar problem for an arbitrary $n$ for the groups of type $B_n$ and $D_n$.\label{nota1}}
As a corollary, we can proof the Isaac's conjecture in this case.
\corop{The number\label{coro:Isaacs_BD} of the coadjoint orbits of maximal dimension of groups of type $U(\Phi)$\textup, where $\Phi=B_n$ or $D_n$ over the finite field $\Fp_q$ is given by the formula $S_n=v^{k_1}(v+1)^{k_2}(v+2)^{\delta}$\textup, where $v=q-1$\textup, and 
$$k_1=\begin{cases}
    (n-1)/2,&\text{ if }\Phi=D_n\text{ or }B_n\text{ and }n\text{ is odd},\\
    (n-2)/2,&\text{ if }\Phi=D_n\text{ and }n\text{ is even},\\
    n/2,&\text{ if }\Phi=B_n\text{ and }n\text{ is even},
\end{cases}$$
$$k_2=\begin{cases}
    (n+1)/2,&\text{ if }\Phi=B_n\text{ and }n\text{ is odd},\\
    (n-2)/2,&\text{ if }\Phi=B_n\text{ and }n\text{ is even},\\
    (n-1)/2,&\text{ if }\Phi=D_n\text{ and }n\text{ is odd},\\
    n/2,&\text{ if }\Phi=D_n\text{ and }n\text{ is even},
\end{cases}$$
for $\delta=1$ in case $\Phi=B_n$ and $n$ is even\textup, and $\delta=0$ otherwise.
}{The proof follows right from the statement of the theorem above by direct calculation.}

\bigskip\textsc{Mikhail Ignatev: National Research University Higher School of Economics,\break Po\-krov\-sky Boulevard 11, 109028, Moscow, Russia}

\emph{E-mail address}: \texttt{mihail.ignatev@gmail.com}

\bigskip\textsc{Mikhail Venchakov: National Research University Higher School of Economics, Pokrovsky Boulevard 11, 109028, Moscow, Russia}

\emph{E-mail address}: \texttt{mihail.venchakov@gmail.com}


\begin{thebibliography}{XXXXX}
\bibitem[Bo03]{Bourbaki03} N. Bourbaki. Algebra II. Chapters 4--7. Elements of mathematics (Berlin). Springer--Verlag, Berlin, 2003.
\bibitem[An95i]{Andre95i} C.A.M. Andre. Basic characters of the
unitriangular group. J. Algebra \textbf{175} (1995), 287--319.
\bibitem[An95ii]{Andre95ii} C.A.M. Andre. Basic sums of coadjoint orbits of the
unitriangular group. J. Algebra \textbf{176} (1995), 959--1000.
\bibitem[AN06]{AndreNeto06} C.A.M. Andre, A.M. Neto. Super-characters of finite
unipotent groups of types $B_n$, $C_n$ and $D_n$. J. Algebra \textbf{305} (2006) 394--429.

\bibitem[GMR16]{GoodwinMoschRoehrle16} S.M. Goodwin, P. Mosch, G. R\"ohrle. On the coadjoint orbits of maximal unipotent subgroups of reductive groups. Transformation Groups \textbf{21} (2016), 399--426.

\bibitem[Hi60]{Higman60} G. Higman. Enumerating $p$-groups. I. Inequalities. Proc. London Math. Soc. (3) \textbf{10} (1960), 24--30.


\bibitem[IP09]{IgnatevPanov09} M.V. Ignatev, A.N. Panov. Coadjoint orbits of the group $\mathrm{UT}(7, K)$. J. Math. Sci. \textbf{156} (2009), no. 2, 292--312.

\bibitem[IP25]{IgnatevPetukhov25} M.V. Ignatev, A.V. Petukhov. Coadjoint orbits of low dimension for nilradicals of Borel subalgebras in classical types, arXiv: \texttt{math.RT/2507.20332}.

\bibitem[Is07]{Isaacs07} I.M. Isaacs. Counting characters of upper triangular groups. J. Algebra \textbf{315} (2007),\break 698--719.

\bibitem[Ka77]{Kazhdan77} D. Kazhdan. Proof of Springer's hypothesis.
Israel J. Math. \textbf{28} (1977), 272--286.

\bibitem[Ki62]{Kirillov62} A.A. Kirillov. Unitary representations of nilpotent Lie groups. Russian Math. Surveys \textbf{17}~(1962), 53--104.

\bibitem[Ki04]{Kirillov04} A.A. Kirillov. Lectures on the orbit method. Grad. Stud. in Math. \textbf{64}, AMS, 2004.


\bibitem[Ko12]{Kostant12} B. Kostant. The cascade of orthogonal roots and the coadjoint structure of the nilradical of a Borel subgroup of a semisimple Lie group. Moscow Math. J. \textbf{12} (2012), no. 3, 605--620.

\bibitem[Ko13]{Kostant13} B. Kostant. Center of $U(\nt)$, cascade of orthogonal roots and a construction of Lipsman--Wolf. In: Lie groups: structure, actions and representations, Progr. in Math. \textbf{306}. Birkh$\ddot{\mathrm{a}}$user, 2013, 163--174.


\bibitem[Le10]{Le10} T. Le. Counting irreducible representations of large degree of the upper triangular groups. J. Algebra \textbf{324} (2010), 1803--1817.

\bibitem[Le74]{Lehrer74} G.I. Lehrer. Discrete series and the unipotent subgroup. Composito
Math. \textbf{28} (1974),\break fasc.~1, 9--19.

\bibitem[Lou11]{Loukaki11} M. Loukaki. Counting characters of small degree in upper unitriangular groups. J. Pure Appl. Algebra \textbf{215} (2011), no. 2, 154--160.


\bibitem[Mb11]{Marberg11} E. Marberg. Combinatorial methods of character enumeration for the unitriangular group. J. Algebra \textbf{345} (2011), 295--323.


\bibitem[Mj97]{Marjoram97} M. Marjoram. Irreducible characters of a Sylow $p$-subgroups of the orthogonal group. Extracta Mathematicae \textbf{12} (1997), no. 3, 315--319.

\bibitem[Mj97']{Marjoram97'} M. Marjoram. Irreducible characters of Sylow p-subgroups of classical groups. PhD thesis,
National University of Ireland, Dublin, 1997.

\bibitem[Mj99]{Marjoram99} M. Marjoram. Irreducible characters of small degree of the unitriangular group. Irish Math. Soc. Bull. \textbf{42} (1999), 21--31.




\bibitem[Ve26]{Venchakov} M. Venchakov. Rook placements and coadjoint orbits for maximal unipotent subgroups of finite symplectic groups, arXiv: \texttt{math.RT/2602.14933}.

\bibitem[VLA03]{VeraLopezArregi03} A. Vera-L\'opez, J.M. Arregi. Conjugacy classes in unitriangular matrices. Linear Algebra and its Applications \textbf{370} (2003), 85--124.



\end{thebibliography}

\begin{thebibliography}{99}

\bibitem{Humpreys1}Хамфри Дж. Введение в теорию алгебр Ли and их
представлений. --- М.: МЦНМО, 2003.
\bibitem{Ignatev}Игнатьев М.В. Субрегулярные характеры
унитреугольной группы над конечным полем. Фундаментальная and
прикладная математика, т. \textbf{13}, вып. \textbf{5}, 2007, с.
103-125, см. также arXiv:math.RT/0801.3079v2.
\bibitem{IgnatevPanov} Игнатьев М.В., Панов А.Н.
Коприсоединённые орбиты группы $\mathrm{UT}(7, K)$. Фундаментальная
and прикладная математика, т. \textbf{13}, вып. \textbf{5}, 2007, с.
127-159, см. также arXiv: math.RT/0603649v3.
\bibitem{Kirillov1} Кириллов А.А. Лекции по методу орбит. ---
Новосибирск: Научная книга ИДМИ, 2002.
\bibitem{Kirillov2} Кириллов А.А. Унитарные представления
нильпотентных групп Ли. УМН, т. \textbf{17}, 1962, с. 57-110.
\bibitem{Kirillov3} Кириллов А.А. Метод орбит and конечные группы. --- М.: МЦНМО, МК НМУ, 1998.
\bibitem{Panov} Панов А.Н. Инволюции в $S_n$ and ассоциированные
коприсоединённые орбиты. Зап. научн. сем. ПОМИ, т. \textbf{349},
2007, с. 150-173, см. также arXiv: math.RT/0801.3022v1.
\end{thebibliography}
\end{document}